\documentclass[12pt]{iopart}
\usepackage{amsfonts}
\usepackage{xfrac}
\usepackage{xcolor}
\usepackage{url}
\usepackage{graphicx,caption}

\newcommand{\R}{\mathbb{R}}

\usepackage[linktocpage=true]{hyperref}

\begin{document}
\title[Knot intensity distribution]{Knot intensity distribution: a local measure of entanglement}

\author{Agnese Barbensi$^1$ and Daniele Celoria$^1$}
\address{$^1$ School of Mathematics and Statistics, University of Melbourne, 3010~VIC, Australia}

\begin{abstract}
The problem of finding robust and effective methods for locating entanglement in embedded curves is relevant to both applications and theoretical investigations. Rather than focusing on an exact determination, we introduce the knot intensity distribution, a local quantifier for the contribution of a curve's region to global entanglement. The integral of the distribution yields a measure of tightness for knots. We compute the distribution for ideal knots, and study its behaviour on prime and composite random knots. Intensity distributions provide an effective method to locate entanglement. In particular, they identify regions in knots that accommodate passages leading to topological changes. 
\end{abstract}

\section{Introduction}
Everyone has an intuitive idea of what it means for a piece of rope to be knotted. A more challenging task is to pinpoint regions in the rope that are actively contributing in creating the knot. This is often relevant in applications; as an example, the folding and stability of knotted proteins is known to be influenced by the geometry and localisation of the entanglement in their native states \cite{dabrowski2016,piejko2020folding,barbensi2021topological}. 

A knotted arc is called \textit{tight} if its entanglement is highly localised, so that the portion involved in the knot is small when compared to the length of the whole curve. To give a precise quantification of tightness, it is necessary to pinpoint the entanglement's location. The most straightforward strategy is to analyse the topology of all the curve's sub-arcs; this procedure determines the shortest knotted portion~\cite{multiscale,knotoid}, which is often called the arc's \textit{knot core}. Tightness can then be measured as a function of length and position of the knot core, relative to the entire curve. 

While being mainly motivated by applications in the open setting, similar ideas generalise to the problem of measuring tightness of closed curves,~\textit{i.e.}~mathematical knots. Standard approaches to this question rely on extending the open arc setting to identify a specific portion of a knot as its knotted core \cite{knotoid, sumners1990detecting,marcone2004length, millett2005linear, tubiana2011probing, rawdon2015subknots}. 

Despite the obvious appeal of looking for an exact identification of the knot core, the concept of the ``shortest knotted portion'' is problematic in general. Indeed, in some cases, the shortest knotted portion might not contain a curve's entire entanglement, and different notions might be better suited to localise it~\cite{multiscale}, see also Figure~\ref{fig:figure1}(A). Here, we explore an alternative point of view; we develop a mathematical framework to quantify the local intensity of the entanglement. To a given a knotted configuration, specified by an embedding $\gamma\colon S^1 \rightarrow \R^3$, we associate a function $f_\gamma\colon S^1 \longrightarrow \R_{\ge 0}$, called the \emph{knot intensity distribution}. This function assigns a value to each point in a curve, measuring the amount of local contribution to the overall topology, computed from a superposition of knot cores. 

The key advantage of this strategy is to remove some intrinsic ambiguities present in previous methodologies, while encoding nuanced information on the topology and geometry of the embedding. This information can be extracted using a number of techniques, including a coarea approach. In this case, the resulting \textit{fingerprint} functions can be used to create a coordinate-free statistical summary of the distributions. As another example, the integral of $f_\gamma$ yields a non-negative number $\mathcal{D}(\gamma)$. We call this number the \emph{knot density} of $\gamma$, and show that it provides a reliable measure of tightness. Indeed, high values of $\mathcal{D}$ arise when the intrinsically knotted component of a knot is evenly distributed throughout the embedding, while small values indicate localised entanglement.  The maximal value of this integral among all representatives of a given knot type $K$ is a knot invariant $\mathcal{D}(K)$; it would be interesting to understand its discriminatory power (\emph{cf}.~the discussion in Section~\ref{sec:results}). 

An explicit computation of $f_\gamma$ and related quantities is not feasible for a generic smooth embedding. We therefore extend the framework to the piece-wise linear (PL) setting, where $f_\gamma$ can be approximated computationally. As a first analysis, we compute the knot intensities distributions for PL representations of \textit{ideal}~\cite{stasiak_nature,ideal_stasiak} embeddings of prime and composite knots with low crossing number. We then examine the behaviour of the knot intensity distributions for random equilateral PL embeddings of low-crossing number knots, as a function of their length. The resulting statistics show a clear stratified behaviour with respect to crossing number, and a tendency of composite knots towards less tight conformations. 

As the distribution quantifies the intensity of knottiness for each sub-arc, it induces a hierarchical classification of points in the curves based on their contribution to global entanglement. We show that this information can be used to discriminate regions based on whether they allocate strand passages leading to topology changes.

\section{Methods}
Our goal is to associate a point-wise measure of entanglement to a knotted embedding. The method we propose is based on considering a closed curve as the superposition of all of its openings, and on identifying a knot core candidate for each of the openings, using standard methods for open curves~\cite{knotoid}. The knot intensity is then computed based on this information. The use of linear chains to characterise the topology of closed knots was also explored in~\cite{rawdon2015subknots} using the so-called \textit{disk matrices}, that encode the entanglement of all of a knot's sub-chains.  

\subsection{Mathematical framework}
The definition of the knot intensity distribution follows three steps. First, \textbf{(1)} we identify the knot type of a smooth embedding of an arc $\alpha\colon [0,1] \longrightarrow \R^3$. Standard methods for this task rely on artificially ``closing'' the curve to transform it into a knot. This \textit{closure} operation can be done in a number of different ways, each coming with strengths and weaknesses~\cite{sumners1990detecting,marcone2004length, millett2005linear, tubiana2011probing, knotoid}. Once such a choice is made, \textbf{(2)} the \emph{knot core} $\mathcal{C}(\alpha)$ of $\alpha$ can be determined by computing the knot type of each sub-arc. Specifically, for every $t \in [0,1]$, we consider the embeddings $\alpha^t\colon [\frac{t}{2}, 1 - \frac{t}{2}] \longrightarrow \R^3$ obtained by restricting the domain of $\alpha$. Then $\mathcal{C}(\alpha)$ is the image of the restriction $\alpha^{t_0}$, where $t_0$ is the minimal value of $t$ such that the closure of $\alpha^t$ has the same associated knot type as $\alpha$ (with respect to the fixed closure choice). Finally, \textbf{(3)} let  $\gamma \colon S^1 \hookrightarrow \R^3$ be a smooth embedding representing a non-trivial knot, with arclength parametrisation.  For $\theta\in S^1$, the \emph{knot intensity distribution} of $\gamma$ at the point $\gamma(\theta)$ is defined as the measure of the set $$\{\phi \in S^1\,|\,\gamma(\theta) \in \mathcal{C}(\gamma(S^1 \setminus (\phi-\epsilon, \phi + \epsilon)))\},$$ for some $\epsilon \ll 0$, renormalised by dividing by $2\pi$. The output of this construction is a function $f_\gamma\colon S^1 \longrightarrow [0,1]$. We note that the choice of a sufficiently small $\epsilon$ does not affect the following discussion. 

\begin{figure}[ht]
\centering
\includegraphics[width = 12cm]{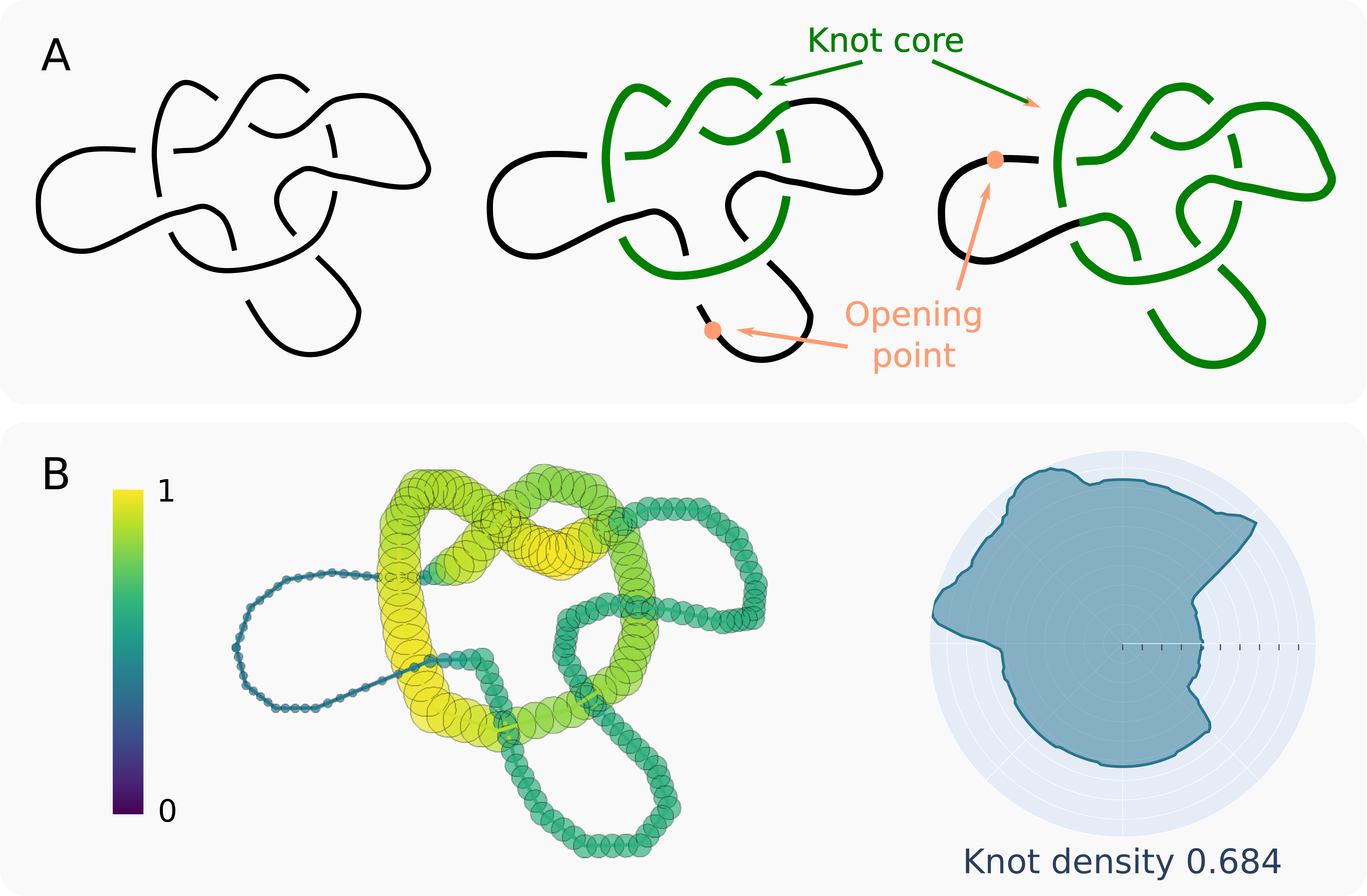}
\caption{\textbf{(A)} On the left, an embedding of the trefoil knot. On the right, the knot cores obtained from two different opening points. \textbf{(B)} On the left, the knot intensity distribution is plotted on the embedding, using colour and thickness. The colour scale will be maintained throughout the paper. On the right, a radar plot of the distribution with the corresponding density value. Note that the three maxima highlight the presence of the embedding's ``essential'' crossings.}
\label{fig:figure1}
\end{figure}

Figure~\ref{fig:figure1}(A) shows an example of an embedding $\gamma$ where different openings result in different knot cores; one corresponds to the shortest arc with the same topology as $\gamma$, while the other to the shortest knotted arc with trivial complement, whose topology doesn't change if extended (\textit{cf.}~\cite{multiscale}). Part~(B) of the same figure conveys instead why our methodology is independent of these issues, while retaining the information contained in previous constructions. Indeed, the two choices of knot core are recovered as a weighted superposition. \\

Given a function $f \colon S^1 \longrightarrow [0,1]$, obtained through the procedure detailed above, define its \emph{fingerprint} as the coarea function $\overline{f} \colon [0,1] \longrightarrow [0,1]$, with $$\overline{f} (t) = \int_{S^1} \min\{f(s),t\} \,\mathrm{d}s.$$ 

In other words, $\overline{f}$ is constructed from the integral of the level sets of $f$. The fingerprint doesn't unique determine $f$, but gives a coordinate-independent way of examining the behaviour of the distributions; fingerprints can be added and averaged, and thus can be used to study the statistical behaviour of the knot intensities we are interested in. The fingerprint has the following elementary properties: $\overline{f}(0) = 0$, $\overline{f}(1) = \int_{S^1} f\, \mathrm{d}s$, $0 \le \overline{f}' \le 1$ and $\overline{f}'' \le 0$. A knot invariant is extracted from the intensity distributions by considering the maximal value of the integral among all embeddings representing a fixed knot type; more precisely, define $$\mathcal{D}(K) = \sup\left\{\int_{S^1} f_\gamma(\theta) \mathrm{d}\theta\,\Big|\,[\gamma(S^1)] = K\right\},$$ where square brackets denote the topological type represented by the curve.

\subsection{The PL setting}

An explicit computation of $f_\gamma$ and all related quantities is not feasible for a generic smooth embedding. We therefore extend the framework to PL curves, as it allows performing computations and statistical analysis. 
We consider equilateral PL curves, that is, closed chains whose segments have the same length. Taking each segment to be of length one, the length of the curve is the number of edges.

For a given PL embedding of length $N$, we consider the $N$ open curves obtained by cyclically removing a single edge from the embedding. The corresponding knot cores are then computed via a probabilistic closure method, combined with a top-down approach to sub-chain analysis. Specifically, we use Knoto-ID's~\cite{knotoid} knot core identification algorithm with probabilistic closures. For each opening, the output consists of the shortest sub-chain obtained by progressively trimming the open PL curve without changing the dominant knot type. Once the knot cores are determined, we compute the distribution as follows: to each vertex in the curve, we associate a value equal to the number of openings containing it in their core. These values are then normalised by dividing by the length.

\section{Results}\label{sec:results}
In our analysis, we consider prime knot types up to $6$ crossings, and non-trivial connected sums up to $7$ crossings. As a first experiment, we compute the knot intensity distributions of PL curves closely approximating ideal embeddings. Then, we analyse randomly generated knotted configurations, and study the statistical behaviour of distribution, density, and fingerprint by knot type and as the length of the configurations increases. Lastly, we show that there is a strong correlation between peaks of the distribution and strand passages disrupting an embedding's topology.

\subsection{Ideal knots}\label{sec:ideal}

Ideal knots are a special type of knot embeddings, arising in the context of biologically inspired constructions~\cite{stasiak_nature,ideal_stasiak}. A knot embedding is \emph{ideal} if it minimises the length over diameter ratio within its isotopy class. It is still unknown whether there is a unique ideal representative for any given knot type; in this context, it is noteworthy to recall that in general, ideal embeddings do not necessarily maximise symmetries~\cite[Ch.~1.5]{ideal_stasiak} nor minimise occupied volume~\cite{klotz}.

\begin{figure}[h!]
\centering
\includegraphics[width=13cm]{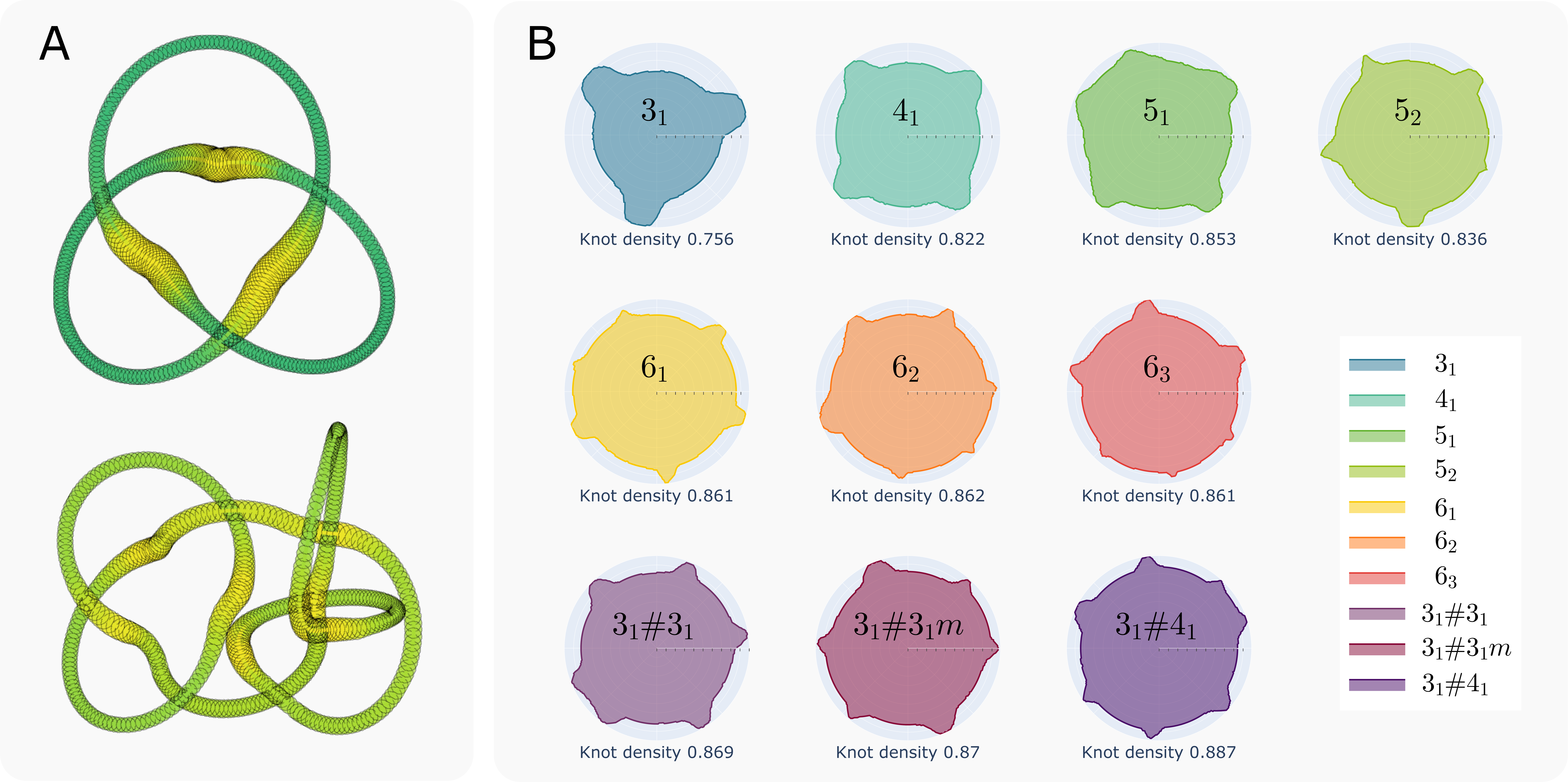}
\caption{\textbf{(A)} Ideal embeddings of the left trefoil $3_1$ (top) and $3_1\# 4_1$ (bottom). Colour and thickness are proportional to knot intensity. \textbf{(B)} Radar plots of the knot intensity distributions for ideal embeddings of prime knots up to $6$ crossings and connected sums (up to mirroring) up to $7$ crossings. The maxima of the distributions correspond to crossings. Values of the knot density are recorded below each plot.}
\label{fig:ideal}
\end{figure}

Equilateral PL embeddings of ideal knots can be generated using the software Ridgerunner~\cite{ridgerunner}; outputs are high-resolution approximations of their smooth counterparts. As a first application of our method, we compute the knot intensity distribution on this dataset. Results are shown in Figure~\ref{fig:ideal}. The value attained by the knot density on ideal knots, shown in Figure~\ref{fig:ideal}(B), is rather high (\textit{cf.}~Section~\ref{sec:random}). This is due to the fact that, by definition, ideal embeddings are optimally using their length to form the knot. That implies that the entanglement is evenly distributed, and even large sub-arcs will likely possess a topology different from the whole embedding.

\subsection{Random knots}\label{sec:random}

Next, we consider a dataset consisting of $10^3$ random configurations, for each of the knot types considered above, and lengths between $50$ and $250$ (with a $50$ step increase). The knotted curves are generated as self-avoiding walks using the Python software Topoly~\cite{topoly}.

\begin{figure}[h!]
\centering
\includegraphics[width = 11cm]{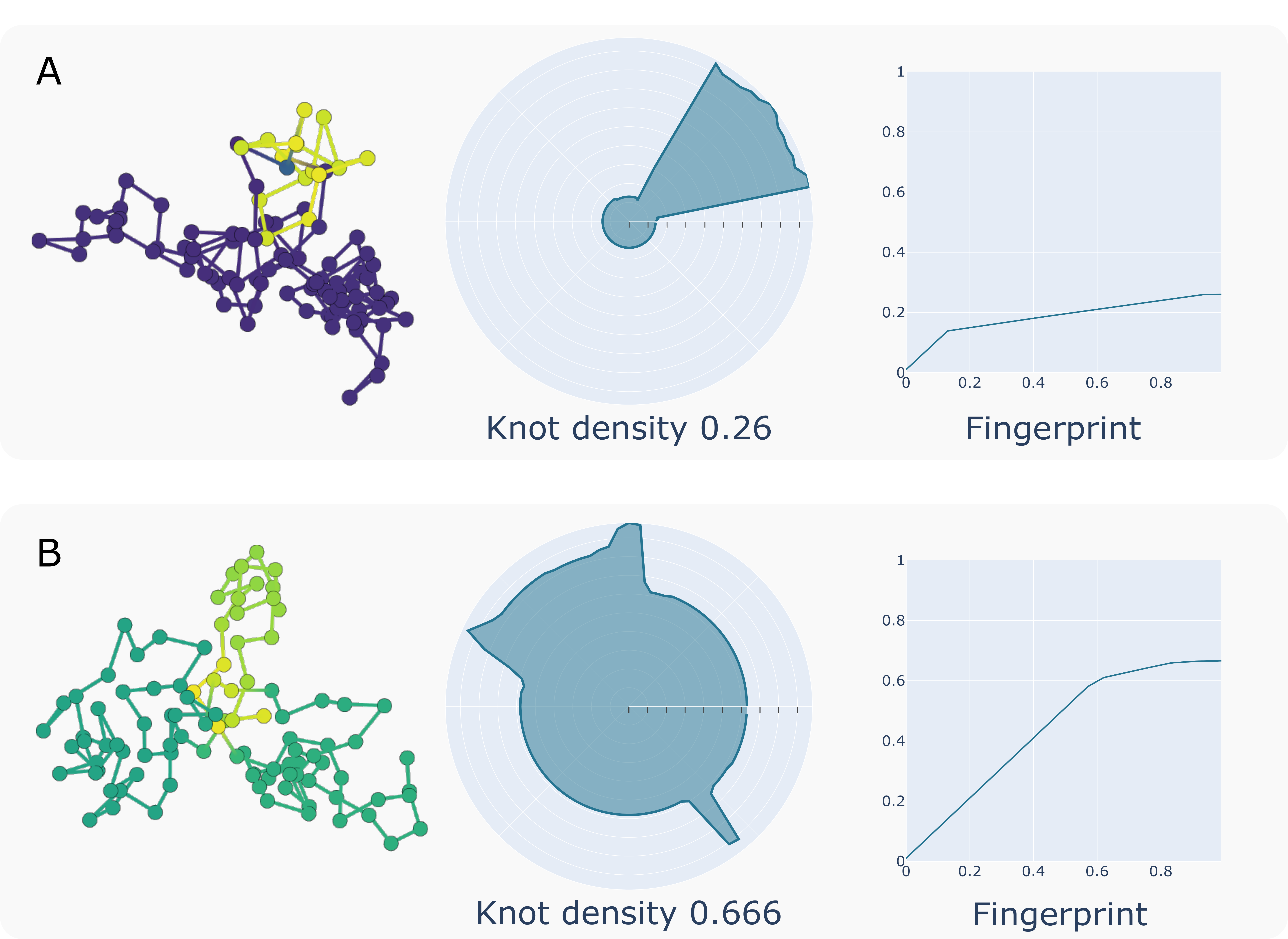}
\caption{\textbf{(A)} A random knotted configuration of length $100$ with $3_1$ topology. Vertices and edges are coloured by the knot intensity. In the middle and right columns, the associated radar plots and knot density values. The entanglement is localised in a small portion of the curve. \textbf{(B)} A random knotted configuration of length $100$ with $3_1$ topology. Vertices and edges are coloured by the knot intensity. In the middle and right columns, the associated radar plots and knot density values. The entanglement is spread throughout a long portion of the embedding.}
\label{fig:trefoil_examples}
\end{figure}

Figure~\ref{fig:trefoil_examples} shows two examples of generated configurations, with trefoil topology and length 100. The configuration in the top of Figure~\ref{fig:trefoil_examples}(A) presents a highly localised entanglement, while the one on the bottom is an example of a shallow trefoil. The knot intensity distributions correctly identify both cases. Indeed, high values of the distribution computed from the curve on the top are concentrated in a short region (see left of (A) and (B)), while the bottom distribution shows a more homogeneous behaviour. Similarly, the density is much higher for the curve on the bottom. In both cases, the fingerprint (Figure~\ref{fig:trefoil_examples}(C)) is linear with slope $1$ until it reaches the minimal value of the distribution.

\begin{figure}[ht]
\begin{center}
\includegraphics[width=14cm]{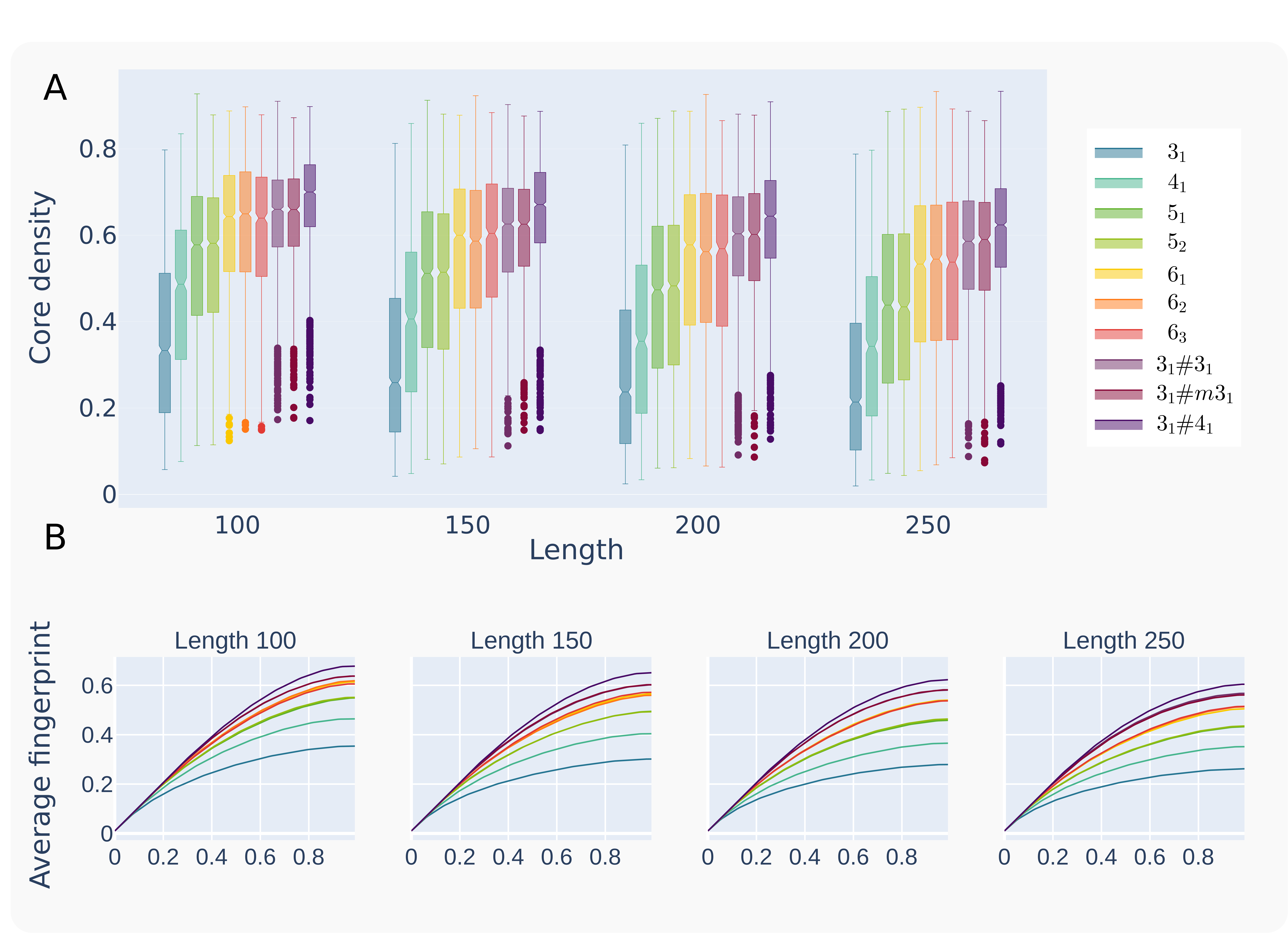}
\caption{The results of the statistical analysis on our random knot population. \textbf{(A)} Box plots of the distributions of knot density, by knot type and length.
\textbf{(B)} Average fingerprints for each knot type. As in (A), there is a clear clustering with respect to crossing number, and separation of the averages densities (\emph{i.e.}~the maxima of the fingerprints). Further, it is apparent from the plots that this stratification persists even when restricting to the (possibly disconnected) sub-chains whose intensity distribution less than a given threshold.}
\label{fig:random}
\end{center}
\end{figure}
To investigate on the behaviour of the knot intensity distribution for different knot types, we subsequently compute the distribution, density and fingerprint of each configuration in the dataset. Results are shown in Figure~\ref{fig:random}. The plot shows a dependency of knot densities on the crossing number. More precisely, prime knots with the same crossing number exhibit very similar knot density distributions, and those distributions are ranked by crossing number, consistently throughout the range of lengths considered. This validates the intuitive idea that more complex knots ``occupy'' longer regions at fixed length. Note that the average density is inversely proportional to the embedding's length; this is consistent with the idea that knot domains tend to be very tight~\cite{zheng2010tightness}. For composite knots, we observe a clear narrowing of the densities with respect to prime knots with the same crossing number, as well as an increase in median values. Again, this is in line with previous studies~\cite{tubiana2014computational} and suggests that prime factors in composite knots tend to spread out rather than merge, leading to a more homogeneous distribution of intensities and a reduced density variation. 

Similarly, the averaged fingerprints shown in Figure~\ref{fig:random}(B) confirm the presence of a clear division by crossing number, which persists at all knot lengths; the same figure implies that the same conclusion holds when restricting to those sub-chains whose intensity distribution values are below a fixed threshold. We also report the absence of correlation between knot density and average writhe, and only a small ($\le$ 0.3) positive correlation with the average crossing number.

The average values of $\mathcal{D}$ on the knots within our population are contained in the plot in Figure~\ref{fig:random}(B), as they correspond to the fingerprints' maxima. The maximal values computed do not present a consistent behaviour for varying length, likely due to noise. It would be interesting to understand to what extent the maximal knot density can distinguish knots; it is reasonable to assume that the values of $\mathcal{D}(K)$ should increase in an approximately monotonic fashion, as a function of the minimal crossing number of $K$. The computations performed in this paper are however not sufficient to support this claim.

We also note that the minimum of the knot density for any given knot type tends to $0$; indeed, it is possible to generate configurations minimising the value of $\mathcal{D}$ within a fixed knot type: namely, consider a PL embedding for a knot $K$ realising its \emph{equilateral stick number} $esn(K)$~\cite{sticknumber}. This quantity is an integer valued invariant of knots, defined to be the minimal number of straight segments of the same length needed to represent a given knot type. Modify this embedding by adding a sequence of $N-est(K)$ further edges in an unknotted fashion; in this modified embedding, the only part contributing to the intensity distribution is entirely concentrated in the portion entailed by the minimal embedding. Therefore, the knot density will be close to $\frac{esn(K)}{N}$, which converges to $0$ for $N \rightarrow \infty$.

\subsection{Knot intensity locates entanglement and non-cosmetic strand passages}\label{sec:crossing_changes}
The knot intensity of low crossing number ideal knots has peaks in bijection with crossings, see Figure~\ref{fig:ideal}(B); a similar behaviour can be observed in Figure~\ref{fig:figure1}(B). A natural question is whether this is a special case of a more general phenomenon. Namely, we hypothesise that higher knot intensity values correlate with regions in the knot whose alteration leads to a change in topology. We test this hypothesis as follows; starting from a random PL trefoil of length $40$, we perform successive random crankshaft moves, allowing non-phantom moves as well. At each strand passage, we record the starting and ending PL curves, and their corresponding knot types. We label as \textit{cosmetic} the strand passages preserving the topology, and discard cosmetic changes between trivial knots. 

\begin{figure}[ht]
\begin{center}
\includegraphics[width=12cm]{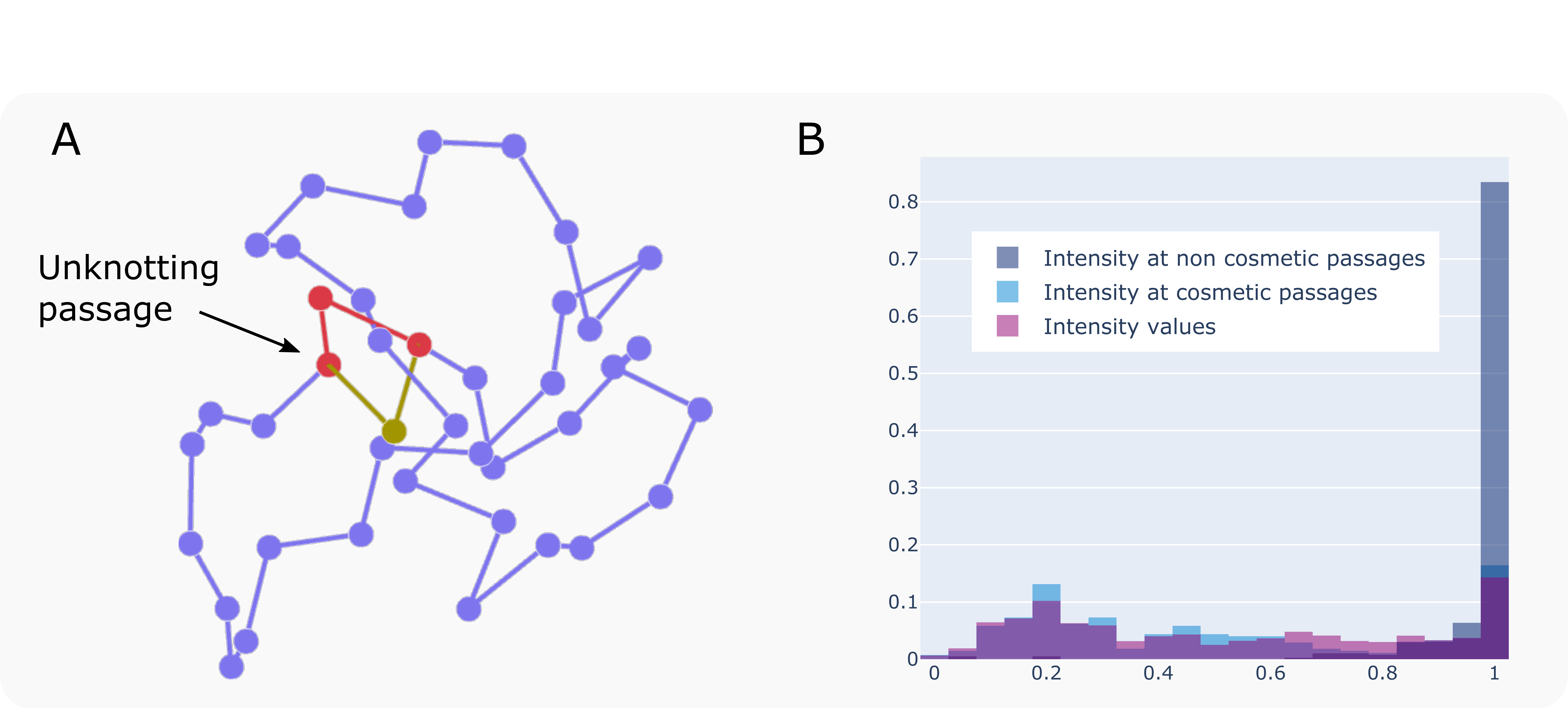}
\caption{\textbf{(A)} The superposition of two PL embeddings, differing by a single non-phantom and non-cosmetic crankshaft move. The two associated knot types are the trefoil and the trivial knot. \textbf{(A)} Bar plots for the normalised values of the intensity distribution (pink), and the same values restricted to crossing change sites in the cosmetic (light blue) and non-cosmetic (dark blue) cases.}
\label{fig:cosmetic}
\end{center}
\end{figure}
Each such pair of embeddings identifies two points, corresponding to the centre of the crankshaft move in each embedding, and to the other segment where the self-intersection occurs. The value of the intensity distribution is then computed at each of these points. We keep the highest among these two values, and normalise it by the maximum of the corresponding distribution. The idea is that if the initial assumption is correct, then the values on non-cosmetic crossing changes should be on average higher than the cosmetic ones. The results of this analysis, displayed in Figure~\ref{fig:cosmetic}, strongly confirm the hypothesis. Indeed, while cosmetic values are distributed similarly to the total distribution, the non-cosmetic ones present a strong bias towards being localised at the maximum. 
In other words, moves leading to passages altering a knot's topology are identified by maximal values of the knot intensity distribution.

\section{Availability and implementation}\label{sec:code}
The GitHub repository \url{https://github.com/agnesebarbensi/knot_intensity_distribution}
contains the code to compute the knot intensity distribution for a given PL knot embedding, as well as the data and code needed to reproduce the results and figures in the paper.

\subsection*{Acknowledgements}
AB gratefully acknowledges funding through a MACSYS Centre Development initiative from the School of Mathematics and Statistics, the Faculty of Science and the Deputy Vice Chancellor Research, University of Melbourne. DC was supported by Hodgson-Rubinstein’s ARC grant DP190102363 ``Classical And Quantum Invariants Of Low-Dimensional Manifolds''.

\section*{References}
\bibliographystyle{unsrt.bst}

\end{document}